\documentstyle[12pt]{article}
\begin{document}
\bibliographystyle{c:/pctex/texbib/plain}
\newtheorem{ct}{\indent Connection}
\newtheorem{th}{\indent  Theorem}
\newtheorem{df}{\indent  Definition}
\newtheorem{prop}{\indent Proposition}
\newtheorem{rk}{\indent Remark}
\newtheorem{lm}{\indent Lemma}
\newtheorem{prb}{\indent Problem}
\newtheorem{con}{\indent Conjecture}
\newtheorem{cy}{\indent Corollary}

\begin{center}{\Large\bf Equations in a free $\bf Q$-group}\\
O. Kharlampovich, A. Myasnikov
\end{center}
\begin{center} {\bf 0. Introduction}\end{center}
Systems of equations over a group have been widely studied
(see, for instance, \cite{ComEd},\cite{GrigKur},\cite{raz}).  This is currently
one of the main streams of combinatorial group theory. 
The problem of deciding if a system of equations in a group has a
solution is a generalization of the word and conjugacy problems.
Makanin \cite{Mak82} 
and
Razborov
\cite{raz} proved one of the most significant results in this area: the 
algorithmic
solvability of systems of equations
in
 free groups. 
Rips and Sela \cite{RipsSela84} solved equations
over hyperbolic groups by reducing the problem to free groups.
Myasnikov and Remeslennikov proved that
 the universal theory is decidable
over free $A$-groups, where $A$ is an  integral domain
of characteristic 0 and  ${\bf Z}$ is a pure subgroup of $A$.

If ${\bf Z}$ is not a pure subgroup of $A$ then the structure of 
a free $A$-group is much more complicated. It turned out (see \cite{BMR}) that 
the crucial case is $A = {\bf Q}$.  Baumslag \cite{baum1} proved that 
the word problem is decidable in free
$\bf Q$-groups. In \cite{MyaKhHyp} we proved that the conjugacy problem in
these groups is decidable.  

A subring $C$ of the ring ${\bf Q}$ is said to be {\em recursive} if 
there is an algorithm which decides whether a given rational number belongs 
to $C$. Any subring of ${\bf Q}$ is of the form ${\bf Q}_{\pi},$ i.e. 
generated by the set $\{\frac{1}{p} \  | p \in \pi \}$, where  $\pi$ is a set 
of primes.  It is not difficult to see that the recursive subrings 
of  ${\bf Q}$ are exactly the rings ${\bf Q}_{\pi}$ for recursive subsets 
$\pi$. If the set $\pi$ is not recursive then the Diophantine problem over 
a free ${\bf Q}_{\pi}$-group $F^{{\bf Q}_{\pi}}$ is undecidable.     
Indeed, let 
$a \in F^{{\bf Q}_{\pi}}$ be an element which ${\bf Q}_{\pi}$-generates its own 
centralizer in $F^{{\bf Q}_{\pi}}$ (i.e. if  $1 \neq a = b^r$, then $r$ is 
invertible in    
${\bf Q}_{\pi}$), then an equation $x^p = a$ has a solution in $F^{{\bf 
Q}_{\pi}}$ 
if 
and only if $p \in \pi$. 

The main result of this paper is the following.
 \begin{th} Let $\pi$ be a recursive set of primes. Then there 
exists an algorithm that decides if a given finite system of 
equations over a free ${\bf Q}_{\pi}$-group
has a solution, and if it does, finds a solution. \end{th}

In particular, the Diophantine problem over a free $\bf Q$-group 
$F^{\bf Q}$ is decidable.

Let $A$ be an arbitrary ring of characteristic 0 with a prime subring ${\bf Z}.$
 The additive isolator 
 $Is_A({\bf Z}) = \{a \in A \mid \exists n(na \in {\bf Z}) \}$
of ${\bf Z}$ in $A$ is a subring of $A$, which is embeddable in ${\bf Q}$. 
Therefore, 
$Is_A({\bf Z}) = {\bf Q}_{\pi(A)}$, where $\pi(A)$ is the set of all prime 
integers, 
which are invertible in $A$.  

Theorem 1 and approximation results from \cite{BMR} imply the following

\begin{th} Let $A$ be a ring   of characteristic 0 
with identity $1$.   Then 
 an algorithm, that decides if a given finite system $W =1$ of 
equations  with coefficients in $F$ 
has a solution in a free $A$-group $F^A$, exists if and only if the set $\pi(A)$
is recursive. Moreover, if $W = 1$ has a solution, the algorithm finds one.

\end{th}

In particular, for any field $K$ of characteristic 0 the Diophantine problem 
with coefficients in $F$ over $F^K$ is decidable.

In order to solve a system of equations in a free ${\bf Q}_{\pi}$-group,
we develop some methods to reduce a system of equations in such a group
to a finite set of systems in a free group. The reduction guarantees that
for every so-called ``minimal'' solution of the system in the free ${\bf  
Q}_{\pi}$-group, 
there exists a solution to at least one of the systems in the free
group satisfying certain conditions restricting the range of values of
the unknowns, and for every such solution to one of the systems in
the free group
there exists a corresponding solution for the original system in  
the free ${\bf Q}_{\pi}$-group. Since by \cite{Mak82}, \cite{Mak84} and 
\cite{raz} it is
possible to decide if a system of equations in a free group has a
solution (satisfying certain restrictions), this reduction allows one to decide 
if 
the system over
the free ${\bf Q}_{\pi}$-group has a solution. This idea
was used for the first time  in \cite{RipsSela84}.  
Notice that a system of equations over a free ${\bf Q}_{\pi}$-group is 
equivalent
to a system of equations (over the  free ${\bf Q}_{\pi}$-group) with 
coefficients
that lie in the free group with the same set of generators. 
Notice also that 
a free ${\bf Q}_{\pi}$-group is a direct limit of hyperbolic groups 
\cite{MyaKhHyp},
but this does not help us directly. Indeed, if a solution of the system
 exists in
the free ${\bf Q}_{\pi}$-group then it belongs to some member of 
the union (direct limit), so 
to some hyperbolic group, but we do not know which 
member.

We can prove the same result as Theorem 1 for a tensor 
${\bf Q}_{\pi}$-completion of an arbitrary torsion-free hyperbolic group,
but a lot of additional work is required. With a view to proving the
more general result in a subsequent paper, we formulate most of the notions
and lemmas in this paper in more general terms than are required for
proving Theorem~1.

Let $G$ be a torsion-free hyperbolic group with generators $d_1,\ldots ,d_N$. 
We will consider a finite system of equations over the 
${\bf Q}_{\pi}$-group $G^{{\bf Q}_{\pi}}$ (see the definition in the next 
section). By adding
a finite number of  new variables and new equations 
we can reduce this system to a system
with coefficients in $G$ (indeed, every constant of the form
$d^{m/n}$ can be replaced by a new variable $z$ satisfying the
equation $z^n=d^m$). We also can replace this system by an equivalent
 system of 
triangular equations (every equation 
contains no more than three terms). We will consider now a finite
system of triangular equations
with coefficients in $G$
\begin{equation}\label{1}
{\phi}_i (d_1,\ldots d_N,x_1,\ldots x_L)=1
\end{equation}

\begin{center} {\bf 1. $A$-groups}\end{center}
Let $A$ be an arbitrary associative ring with identity and $G$ a group.
Fix an action of the ring $A$ on $G$, i.e. a map $G \times A \rightarrow
G$. The result of the action of $\alpha \in A$ on $g \in G$ is written
as $g^\alpha$. Consider the following axioms:

\begin{enumerate}
\item $g^1=g$, $g^0=1$,  $1^\alpha = 1$ ;
\item $g^{\alpha +\beta}=g^\alpha \cdot g^\beta, \ g^{\alpha \beta} =
(g^\alpha)^\beta$;
\item $(h^{-1}gh)^\alpha = h^{-1}g^\alpha h;$
\item $[g,h]=1 \Longrightarrow (gh)^\alpha = g^\alpha h^\alpha.$
\end{enumerate}

\begin{df}
Groups with $A$-actions satisfying axioms 1)--4) are called {\em $A$--groups}. 
\end{df}
In particular, an arbitrary group $G$ is a ${\bf Z}$-group. We now recall
the definition of $A$-completion. 

\begin{df}

Let $G$ be a group . 
Then an ${A}$--group $G^{A}$ together with a homomorphism 
$G\rightarrow G^A$ is called  a {\em tensor ${A}$}--completion
of the group $G$ if $G^{A}$ satisfies the following universal
property:
for any ${A}$--group $H$ and a homomorphism 
$\varphi: G \rightarrow H$
there exists a unique ${A}$--homomorphism $\psi: G^{A} \rightarrow H$
(a homomorphism that commutes with the action of ${A}$) 
such that the following diagram commutes:

\medskip

\begin{center}

\begin{picture}(100,100)(0,0)
\put(0,100){$G$}
\put(100,100){$G^{A}$}
\put(0,0){$H$}
\put(15,103){\vector(1,0){80}}
\put(5,93){\vector(0,-1){78}}
\put(95,95){\vector(-1,-1){80}}
\put(-10,50){$\phi$}
\put(55,45){$\psi$}
\put(50,108){$\lambda$}
\end{picture}

\end{center}

\medskip
\end{df}
It was proved in \cite{bau} 
that for every group $G$  the tensor ${\bf Q}_{\pi}$-completion of $G$ 
exists and is unique; in \cite{MyasExpo2} this fact was proved for an arbitrary 
ring 
$A.$  

 We describe below the ${\bf Q}_{\pi}$-completion $G^{{\bf Q}_{\pi}}$ 
of a torsion-free
hyperbolic group $G$ as the union of an effective chain of hyperbolic subgroups
(details can be found in \cite{MyaKhHyp}).

An  element $v$ of a group is called a 
{\em primitive} element if it is not a proper power. 

Let ${\bf Z}_{\pi}$ be the multiplicative submonoid of $\bf Z$
generated by the set $\pi $, ${\bf Z}_{\pi}=\{m_1,m_2,\ldots \},$ where 
$m_1<m_2<\ldots .$
For an arbitrary torsion-free hyperbolic group $R$ and natural number $n\geq 2$ 
choose a set of elements
${\cal V}_n = \{v_1\ldots v_t \}\in R$ satisfying the following condition
${\em (S_n)}$:
\begin{enumerate}
\item[1)] ${\cal V}_n$ consists of cyclically minimal (of minimal
length in its conjugacy class) 
primitive elements of length not more than $m_n$;
\item[2)] no two centralizers in the set of centralizers $
\{C(v), v\in {\cal V}_n\}$ are conjugate in $R;$
\item[3)] the set ${\cal V}_n$ is maximal with properties 1) and 2), i.e.
any element of length not more then $m_n$ is conjugate to a power of some
 $v\in {\cal V}_n$.
\end{enumerate}
By definition, 
\begin{equation}
\label{(2)}R({\cal V}_n)=(\ldots 
(R\ast_{v_1=w_1^{m_n}}<w_1>)\ast_{v_2=w_2^{m_n}}<w_2>)\ast
\ldots )\ast_{v_t=w_t^{m_n}}<w_t>). \end{equation}
Notice that this definition does not depend on the order of elements in
${\cal V}_n.$

It was proved in \cite{MyaKhHyp} that $G^{{\bf Q}_{\pi}}$ is the union of a 
chain 
of 
hyperbolic groups 
$$G=T_{0}<T_1<T_2<\ldots < \bigcup _{n=0}^{\infty}T_n,$$
with $T_n=T_{n-1}({\cal V}_n),$ where ${\cal V}_n$ satisfies the
condition $S_n$ in the group $T_{n-1}.$  

\begin{df} Let a group $H$ be an amalgamated product $H=S\ast _{w=t^r}<t>$, 
then
 {\em $t$-syllables} of the word 
$b_0t^{{\alpha} _1}b_1t^{{\alpha} _2}\ldots
b_n,$ where $b_j\in S,$ are the subwords $  
t^{{\alpha} _1},\ldots  ,t^{{\alpha} _n}$. If $u$ is an element in
$H$, then
$|u|_{H}$ is
the number of occurences of $<t>$-syllables  in a reduced word
representing $u$. We call this number
the {\em $t$-length} of $u$.\end{df}

Any finite set of elements $\{g_1,\ldots ,g_n\}$ in $G^{{\bf Q}_{\pi}}$   
is contained in some subgroup $H$ that is obtained from the group $G$
by adding a finite number of roots.     It is the
union of a chain of subgroups $H_i$. The groups $H_i$ together with
a canonical set of generators are defined below. Let $G=H_0.$
\begin{enumerate}
\item Step 1. Consider pairwise nonconjugated cyclically minimal 
primitive elements $u_1,\ldots ,u_{k_1}\in G,$ $|u_1|\leq \ldots \leq |u_{k_1}|$
(here $|u|$ denote the length of $u$ in $G$), 
and add roots
$t_1,\ldots ,t_{k_1},$ such that $u_j=t_j^{s_j}.$ (Notice that $u_{i+1}$
does not become a proper power after we add roots $t_1,\ldots ,t_i$.)
The corresponding groups are
denoted 
by $H_1,\ldots ,H_{k_1},$ where $H_{j+1}=H_j\ast _{u_{j+1}=
t_{j+1}^{s_{j+1}}}<t_{j+1}>.$  
\item Step 2. Consider pairwise nonconjugated primitive elements 
$u_{k_1+1},\ldots ,u_{k_2}\in H_1,$ 
cyclically
reduced in the amalgamated product, each having the reduced 
form $u=t_1^{{\alpha}_1}c_1\ldots t_1^{{\alpha}_k}c_k,$ where $\alpha
_i\not = 0, \alpha _i\in {\bf Z},
 c_i\in G,$
 $|u_{k_1+1}|_{H_1}\leq\ldots \leq |u_{k_2}|_{H_1};$ 
 and  add roots
$t_{k_1+1},\ldots ,t_{k_2},$ 
such that $u_j=t_j^{s_j},$ to the group 
$H_{k_1}$. The corresponding groups are denoted 
by $H_{k_1+1},\ldots ,H_{k_2},$ where $H_{j+1}=H_j\ast _{u_{j+1}=
t_{j+1}^{s_{j+1}}}<t_{j+1}>.$  
\item Step $i+1$. Suppose that $H_1,\ldots, H_{k_i}$ have been constructed.

Consider pairwise nonconjugated primitive 
elements $u_{k_i+1},\ldots ,u_{k_{i+1}}\in H_i$ 
, cyclically
reduced in the amalgamated product, each having the following reduced 
form $t_i^{{\alpha}_1}c_1\ldots t_i^{{\alpha}_k}c_k,$ 
where $\alpha _i\not = 0, \alpha _i\in {\bf Z},
c_i\in H_{i-1}$ ( $c_k$ is not
a power of $u_i,$ because the elements are cyclically reduced), 
$|u_{k_i+1}|_{H_i}\leq\ldots \leq |u_{k_{i+1}}|_{H_i}$ 
,  and  add roots
$t_{k_i+1},\ldots ,t_{k_{i+1}},$ 
such that $u_j=t_j^{s_j},$ to the group 
$H_{k_i}$. The corresponding groups are denoted 
by $H_{k_i+1},\ldots ,H_{k_{i+1}},$ where  
$H_{j+1}=H_j\ast _{u_{j+1}=
t_{j+1}^{s_{j+1}}}<t_{j+1}>.$
\end{enumerate}
Finally, for some number $i$ one has $H=H_{k_{i+1}}.$

The canonical set of generators of $H_0=G$ is $\{d_1,\ldots , d_N\}$,
the canonical set of generators of $H_{j+1}$ is defined inductively
as the union of the canonical set of generators of $H_j$ and ${t_{j+1}}.$

The group $H_{k_{i+1}}$ is called the {\em group at level $i,$}
corresponding to the sequence $u_1,\ldots ,u_{k_{i+1}}.$
The group $H_i$ will be called the {\em group of rank $i.$}
We also order the set of $t_j$'s: $t_k<t_l$ if $k<l.$

Let $F$ be the free group with the same set of generators as $G$
and $F_{k_j+p}=F*K_{k_j+p}$, where $j\leq i$ and $K_{k_j+p}$
is the free group with the generators $t_1,\ldots ,t_{k_j+p}$.
Let
$\pi $ be the natural
homomorphism of $F_{k_i+p}$ onto $H_{k_i+p}.$ Is $v\in F_{j}$ then
by $\tilde v$ we denote $\pi (v)\in H_j.$
\begin{center}{\bf 2. Diagrams}\end{center}
Recall that a {\em map} is a finite,
planar connected $2$-complex. 

By a {\em diagram $\Delta$ over a presentation 
$<a_1,\ldots ,a_m |
R_1,\ldots ,R_n>,$} where the words $R_i$ are cyclically reduced,
we mean a map  with a function $\phi$ 
which assigns to each edge of the map one of the letters $a_k^{\pm 1},$
$1\leq k\leq m,$ such that $\phi (e^{-1})= ({\phi (e)})^{-1}$ and if
$p=e_1\ldots e_d$ is the contour of some cell $\Phi$ of $\Delta,$
then $\phi (p)=\phi (e_1)\ldots \phi (e_d)$ in the free group $F(a_1,\ldots 
,a_m)$ 
is a cyclic
shift of one of the defining words $R_i^{\pm 1}.$ In general the word
$\phi (p)$ is called the {\em label} of the path $p.$ The label of a diagram
$\Delta$ (whose contour is always taken with a counterclockwise
orientation) is defined analogously.

Van Kampen's Lemma states that a word $W$ represents the identity
of the group $G$ if and only if there is a simply   connected
(or Van-Kampen, or disk) diagram $\Delta$ over $G$ such that the boundary
label of $\Delta$ is $W.$

Suppose we have a diagram over $H_j.$ A {\em $t_j$-strip} is a subdiagram with
the
boundary label $t_j^{s_jn}u_j^{-n}$ (see Fig.~1a),
consisting of cells with the boundary $t_j^{s_j}u_j^{-1}$ (see Fig. 1b).
Two $t_j$-strips can be glued together to form a  {\em
paired $t_j$-strip} (see Fig. 1c).
More than two  $t_j$-strips can be glued together to form a {\em
$t_j$-star} (see Fig. 1d).

Every minimal (with minimal number of $t_j$-cells) 
diagram over $H_j$ consists
of paired $t_j$-strips, $t_j$-stars, $t_j$-strips on the boundary 
and $H_{j-1}$-subdiagrams between them
(annular paired $t_j$-strips can be assumed not to occur).
Suppose we have a $H_{k_i}$-diagram. 
 Notice that for $i\leq j,l\leq k_i,$
$t_j$-stars, $t_l$-stars, and paired strips cannot meet as in Fig. 1e, 
because
$u_i,\ldots ,u_{k_i}$ do not contain $t_i,\ldots ,t_{k_i}$.

For every word $w$ which represents the identity element in $H_{k_i}$ there is
a  
diagram over $H_{k_i}$ with the boundary label $w$ that has the form 
shown in Fig. 1f. It consists
of glued $t_j$-strips, and $t_j$-strips on the boundary 
for $j\in\{i,\ldots ,k_i\},$ 
and $H_{i-1}$-subdiagrams between them.
\begin{center}{\bf 3. Some properties of the Cayley graph of $H.$}
\end{center}

A generating set $J$ of $H=H_{k_{i+1}}$ consists of $d_1,\ldots ,d_N$ and 
the added roots $t_1,\ldots , t_{k_{i+1}}$. 
Recall that the vertices of the Cayley graph 
$\Gamma ({H})=\Gamma ({H},J)$ are elements of ${H}$; and
two vertices $g,\ h,$ are connected by an edge $e=(g,d),$ 
with
label
$\phi (e)=d\in J,$ if $h=gd.$
\begin{df} Let $u$ be a cyclically reduced word in $H_i$.
A word $X$ is called a $u$-periodic word if it is a subword of some power
$u^k.$\end{df}
\begin{lm}\label{nar}
Let $p,s\in H_{i-1}$ ,   $u,v\in F_i$ be words representing the
reduced forms of elements $\tilde u,\tilde v\in
\{u_{k_i+1},\ldots ,u_{k_{i+1}}\}$,  
let $X,Y$ be $u-$ and $v-$periodic words respectively and $s\tilde
Xp=\tilde Y.$ 

If $\tilde u\not = \tilde v,$
then
$|{\tilde X}|_{H_i},|\tilde Y|_{H_i}
< |\tilde u|_{H_i}+|\tilde v|_{H_i}+2$.

If  
$\tilde u=\tilde v$ and $X,Y$  begin and end with $t_i$ 
belonging to distinct
$<t_i>$-syllables of $u$, 
$s\not =1$, 
then  $|\tilde X|_{H_i},|\tilde Y|_{H_i}
\leq |\tilde u|_{H_i}$.   

 If  $\tilde u=\tilde v$,
$X,Y$  begin and end with $t_i$
belonging to the same
$<t_i>$-syllables of $u$, 
$s\not =1,$ then $s=u_i^{\beta}, p=u_i^{\gamma}$ in $H_i$ and
the canonical image of $t_i^{\beta s_i}Xt_i^{\gamma s_i}$ equals  
the canonical image of $Y$ in $H_{i-1}*<t_i>.$
\end{lm}
{\bf Proof.} We will prove the first assertion of the lemma.
Suppose $|\tilde v|_{H_i}\geq |\tilde u|_{H_i}.$ 
Taking a cyclic permutation of $u$ instead of $u$
and a cyclic permutation of $v$ instead of $v$ we can suppose that
$X$ starts with $u$, $Y$ starts with $v$, 
$\tilde u=a_0{t_i}^{{\beta}_1}a_1\ldots {t_i}^{{\beta}_n}a_n$, $0<{\beta}_k<s_i$ 
for all $k,$ 
and
$\tilde v=s(a_0{t_i}^{{\beta}_1}a_1\ldots {t_i}^{{\beta}_n}a_n)^ca_0\ldots 
{t_i}^{{\bar\beta}_m}b.$

Suppose that
$|\tilde X|_{H_i}
\geq |\tilde u|_{H_i}+|\tilde v|_{H_i}+2. $
The case $\bar\beta _m<\beta _m$ is impossible, because it implies 
$\beta _1=\beta _m=\beta _m-\bar\beta _m$.
We have
\begin{equation}
\label{eq2} 
bsa_0{t_i}^{{\beta}_1}a_1\ldots {t_i}^{{\beta}_n}a_na_0{t_i}^{{\beta}_1}a_1
{t_i}^{{\beta}_2}
\ldots =a_{m}{t_i}^{{\beta}_{m+1}}a_{m+1}\ldots 
{t_i}^{{\beta}_n}a_na_0{t_i}^{{\beta}_1}a_1
\ldots a_{m+1}{t_i}^{{\beta}_{m+2}}\ldots \end{equation}
Now instead of $a_r, r<n, $ we write $\bar a_r$ and instead
of $a_0a_n$ we just write $\bar a_n.$ 
We have ${\beta}_j={\beta}_{j+m}$ for any $j$ 
(indices are taken modulo $n$).
These equations for powers imply $a_0\not = 1$ or $a_n\not = 1.$
Let $d=(n,m)$ then this implies 
that ${\beta}_j={\beta}_{j+d}$ (indices are taken modulo $n$).
We also have from Equation (\ref{eq2})
$bsa_0=a_mu_i^{\alpha _0},$
$u_i^{\alpha _{j-1}}\bar a_{j}=\bar a_{m+j}u_i^{\alpha _j}, j=1,\ldots ,n$
(the subscripts of the $\bar a$'s are taken modulo $n$).
Then $\alpha _n=\alpha _0,$ because
 $u_i^{\alpha _{0}}\bar a_{1}=\bar a_{m+1}u_i^{\alpha _1}$ and
 $u_i^{\alpha _{n}}\bar a_{1}=\bar a_{m+1}u_i^{\alpha _{n+1}}$
and the subgroup $<u_i>$ is malnormal.
Hence $$u_i^{\alpha _0}a_1\ldots {t_i}^{{\beta}_n}a_na_0{
t_i}^{{\beta}_1}u_i^{-{\alpha _0}}=
a_{m+1}\ldots {t_i}^{{\beta}_n}a_na_0{t_i}^{{\beta}_1}a_1
\ldots a_{m}{t_i}^{{\beta}_{m+1}}.$$
We have for some $\beta$ and $\gamma ,$
$u_i^{\beta}\bar a_d=\bar a_nu_i^{\gamma}.$
We also have $\Sigma_k {\alpha}_{j+dk}=0.$ And $\tilde u$ is the $(n/d)$-th 
power
of the element $a_0t_i^{{\beta}_1}\ldots t_i^{{\beta}_{d}}u_i^{-{\beta}}a_n
=a_0t_i^{{\beta}_1}\ldots t_i^{{\beta}_{d}}a_du_i^{-{\gamma}}.$ 
Then $\tilde v$ is also a proper power.

In the case $m=0, |\tilde v|_{H_i}>|\tilde u|_{H_i},$ we have $\tilde v=
s(a_0{t_i}^{{\beta}_1}a_1\ldots {t_i}^{{\beta}_n}a_n)^cb$
and 
\begin{equation}
\label{eq2'} 
bsa_0{t_i}^{{\beta}_1}a_1\ldots {t_i}^{{\beta}_n}a_na_0{t_i}^{{\beta}_1}a_1
{t_i}^{{\beta}_2}
\ldots =a_{0}{t_i}^{{\beta}_{1}}a_{1}\ldots {t_i}^{{\beta}_n}
a_na_0{t_i}^{{\beta}_1}a_1
{t_i}^{{\beta}_{2}}\ldots \end{equation}
Then $bsa_0=a_0u_i^{\alpha _0}$ and 
$u_i^{\alpha_{j-1}}\bar a_j= \bar a_ju_i^{\alpha_j}.$
This implies that $\alpha _0=0$ and $b=s^{-1}$ and
$\tilde v$ is conjugated to a power of $\tilde u.$ 

In the case $|\tilde v|_{H_i}=|\tilde u|_{H_i},$ we have
$\tilde v=
sa_0{t_i}^{{\beta}_1}a_1\ldots {t_i}^{{\beta}_n}a_nb$
and $$
bsa_0{t_i}^{{\beta}_1}a_1\ldots {t_i}^{{\beta}_n}a_nbsa_0{t_i}^{{\beta}_1}a_1
{t_i}^{{\beta}_2}
\ldots =a_{0}{t_i}^{{\beta}_{1}}a_{1}\ldots {t_i}^{{\beta}_n}a_na_0
{t_i}^{{\beta}_1}a_1
{t_i}^{{\beta}_{2}}\ldots $$
and again $b=s^{-1}$ and $\tilde u$ and $\tilde v$ are conjugated.

The second and third assertions of the lemma can be proved similarly.

\begin{cy}
Let $p,s\in H_{i-1}$  and $\tilde u,\tilde v\in
\{u_{k_i+1},\ldots ,u_{k_{i+1}}\}$, $\tilde u\not = \tilde v$, 
$|\tilde v|_{H_i}\geq |\tilde u|_{H_i}.$ 
Let $X,Y$ be $u-$ and $v-$periodic words respectively 
and $sXp=Y,$ then  $|\tilde X|_{H_i},|\tilde Y|_{H_i}<3|\tilde v|_{H_i}.$
\end{cy}

This follows directly from the lemma in  the case $|\tilde v|_{H_i}>1$. If
$|\tilde v|_{H_i}=|\tilde u|_{H_i}=1$, then taking cyclic permutations of $u,v$ 
instead of $u$ and $v$,
we have $\tilde u=a_0t^{\beta}a_1, \tilde v=sa_0t^{\beta}b_1.$
If $|\tilde X|_{H_i}=|\tilde Y_{H_i}|\geq 3,$ then
$$sa_0t_i^{\beta}a_1a_0t_i^{\beta}a_1a_0t_i^{\beta}\ldots =
sa_0t_i^{\beta}b_1sa_0t_i^{\beta}b_1sa_0t_i^{\beta}\ldots .$$
Hence $a_1a_0=b_1sa_0u_i^{\alpha _0}, 
u_i^{\alpha _0}a_1a_0=b_1sa_0u_i^{\alpha _1}$ for some $\alpha _0,
\alpha _1\in {\bf Z}.$ Then $\alpha _0=\alpha _1=0$, $a_1a_0=b_1sa_0$
and $\tilde u=\tilde v,$ a contradiction.

\begin{df} 
An element in $H_{i+r}$ is said to be written in  {\em reduced
form in rank $i$} if it belongs to
$H_i$ and is in reduced form as an element in the amalgamated product 
$H_{i}=H_{i-1}\ast _{u_{i}=
t_{i}^{s_{i}}}<t_{i}>.$ If $r=0,$ then the element is said to be in
reduced form in all ranks $\geq i$ if it is in reduced form in rank $i$.
An element $h$ in $H_{i+r}$ is defined by induction on $r$ to be written in  
{\em 
reduced
form in all ranks $\geq i$
} if it is written in the reduced form in the amalgamated product
$H_{i+r}=H_{i+r-1}\ast _{u_{i+r}=
t_{i+r}^{s_{i+r}}}<t_{i+r}>,$
$\ \ h= b_0t_{i+r}^{{\alpha} _1}b_1t_{i+r}^{{\alpha} _2}\ldots
b_n$, where the $b_1,\ldots ,b_n\in H_{i+r-1}$  are in the reduced
form in 
all ranks 
$\geq i$. \end{df}

\begin{df} A {\em $<t_i>$-syllable of a path} is a subpath labelled by some 
$<t_i>$-syllable of the label of the path.\end{df}

\begin{df}
A {\em $u$-path} is a path labelled by a $u$-periodic word.\end{df} 
\begin{df}\label{min}
Let $\tilde u, \tilde v\in\{u_{k_i+1},\ldots ,u_{k_{i+1}}\}, |\tilde 
u|_{H_i}\geq 
|\tilde v|_{H_i}.$ 
Consider in $\Gamma (H)$ two paths: a $u$-path $r_1$ and a $v$-path $r_2$, where
$r_1$ connects the  sequence of vertices
$$\ldots , g\tilde u^{-2}, g\tilde u^{-1}, g, g\tilde u, g\tilde u^2,\ldots ,$$
 $r_2$ connects the sequence of vertices 
$$\ldots , gh\tilde v^{-2}, gh\tilde v^{-1}, gh, gh\tilde v, gh\tilde v^2,\ldots 
,$$
the label of every subpath of $r_1$ or $r_2$ is a reduced word in $H_i,$
and $h\in H_{k_i}.$ 
Then a path $q$ that connects $r_1$ with $r_2$ is called a {\em minimal path}
if 
\begin{enumerate}
\item $\phi (q)$ represents an element in  reduced form in all ranks 
greater than or equal to $i$,
\item the number of $<t_i>$-syllables
in $\phi (q)$ 
is minimal for all paths connecting $r_1$ with $r_2.$
\end{enumerate}
\end{df}
\begin{lm}\label{Cayley}
Let $r_1$ and $r_2$ be the paths from Definition \ref{min}, $h\in
H_{k_i}$. 
Then there are two possibilities:
\begin{enumerate}
\item $h$ contains some $t_j$ that is greater than or equal to $t_i$.
Then  there are a uniquely determined number $s_1$  and consecutive 
$<t_i>$-syllables
$z_1, z_2$ and $z_3$ of the subpath of $r_1,$ with the label $u,$ 
between $g\tilde u^{s_1}$
and $g\tilde u^{s_1+2}$, 
such that the initial points of all the minimal paths connecting $r_1$
with $r_2$ 
belong to the subpath of $r_1$ joining $z_1$ with $z_3.$ (And, similarly,
there are a uniquely determined number $s_2$  and consequtive $<t_i>$-syllables
$z_4, z_5$ and $z_6$ of the subpath of $r_2,$ with the label $v,$ 
between $gh\tilde v^{s_2}$
and $gh\tilde v^{s_2+2}$, 
such that the terminal points of all the minimal paths connecting $r_1$
with $r_2$ belong to the subpath of $r_2$ joining $z_4$ with $z_6.$)

\item  $h$ does not contain any $t_j$ that is greater than or equal to $t_i.$ 
Then either the same conclusion is true as in the previous case,
or, for any two paths
$q$ and $q'$ connecting $r_1$ with $r_2$ and such that $|q|_{H_i}=
|q'|_{H_i}=0,$
  we have the equality $\pi\phi ({q^{-1}})\tilde X\pi\phi ({q'})=\tilde Y$, 
where
$X$ is a $u$- and $Y$ is a $v$-periodic word corresponding to the subpaths of 
$r_1$ and
$r_2$ connecting the initial and terminal points of $q$ and $q'$ respectively
. Then $|\tilde X|_{H_i},|\tilde Y|_{H_i}\leq |\tilde u|_{H_i}+|\tilde 
v|_{H_i}+1.$
\end{enumerate}

The numbers $s_1,s_2$ above do not depend on $g.$
\end{lm}

{\bf Proof.} 
The assertion of the lemma in the second case follows 
directly from Lemma~\ref{nar}. 

To prove the assertion in the first case, suppose that 
$q$ and $\bar q$ are two minimal paths connecting $r_1$ with $r_2.$
Let  $p_1$, $\bar p_1$ be their initial points and  $p_2$, $\bar p_2$ 
be their terminal points.
We have to prove that the each of the paths $p_1\bar p_1$ 
and $p_2\bar p_2$ contains at most one
$<t_i>$-syllable. 
 
Consider a minimal  diagram $\Theta ,$ 
with contour
consisting of the subpath of $r_1$ between the points $p_1$ and $\bar p_1$,
the path $\bar q$, the subpath of $r_2^{-1}$ between $\bar p_2$ and $p_2$,
and the path $q^{-1}$. 

Every subpath of $r_1$ or $r_2$ is reduced (i.e. has a reduced label), and the
 elements $q$ and 
$\bar q$ are reduced, hence there cannot be $t_j$-arcs and $t_j$-stars
with two ends on the same side of $\Theta $, as in Fig. 2a,
for $j\geq i$. The diagram $\Theta$ has the form shown in Fig. 2b,
where the strips are paired $t_j$-strips for $j\geq i$ , the 3-ended stars
are $t_i$-stars and the regions  
between
the strips and stars must be $H_{i-1}$ subdiagrams. 
Hence the paths $p_1\bar p_1$ and
$p_2\bar p_2$ cannot contain more than one $t_i$-syllable each.

The lemma is proven.

The lemma immediately implies the following result.

\begin{cy}
\label{Cayley1}
Let $\tilde u=u_{k_i+p}, h\in H_{k_{i}+p-1}.$ 
Consider in $\Gamma (H)$ two $u$-paths $r_1$ and $r_2$,
where $r_1$ connects the  sequence of vertices
$$\ldots , g\tilde u^{-2}, g\tilde u^{-1}, g, g\tilde u, g\tilde u^2,\ldots ,$$
$r_2$ connects the sequence of vertices 
$$\ldots , gh\tilde u^{-2}, gh\tilde u^{-1}, gh, gh\tilde u, gh\tilde
u^2,\ldots ,$$
and every subpath of $r_1$ or  $r_2$ is labelled by a reduced word in $H_i.$

Then there are uniquely determined numbers $s_1, s_2,$ 
such that for every path $\bar q$ connecting some $g\tilde u^{s_3}$
with $gh\tilde u^{s_4}$
\begin{enumerate}
\item $\phi (\bar q)=\phi (q_1)\phi (q)\phi (q_2),$ where $q_1$ is a $u$-path
connecting
$g\tilde u^{s_3}$ with some vertex $p_1$ on  the subpath of $r_1$
joining  $g\tilde u^{s_1}$
  and $g\tilde u^{s_1+3}$ 
(see Fig. 3.),
$q$ is a path connecting $p_1$ with a vertex $p_2$ 
 on the subpath of
$r_2$ joining
 $gh\tilde u^{s_2}$  and $gh\tilde u^{s_2+3},$
$q_2$ is a $u$-path connecting $p_2$ with $gh\tilde u^{s_4},$ \item $q$ 
is reduced in all ranks $j\geq i$
\item the sum of $t_i$-lengths
of pieces of $q$ between $<t_j>$-syllables, for $j>i,$
is minimal for all paths $q$ with the above property.
\end{enumerate}
These numbers $s_1,s_2$ and the labels of the paths $q_1,q,q_2$ depend only on 
$h$
and $u$, and  not  on $g.$

There is a uniquely determined number $s_1,$
such that for every path $\bar q$ connecting some vertex $g$ with $ghu^{s_3},$ 
\begin{enumerate}\item 
$\phi (\bar q)=\phi (q_3)\phi (q_4),$ where $q_4$ connects
$g$ with some vertex $p_1$  on the subpath of $r_2$ 
joining $g\tilde u^{s_1}$  and $g\tilde u^{s_1+3}$
,
$q_4 $ is a $u$-path  connecting $p_1$ with $gh\tilde u^{s_3},$ 
\item $q_3$ is reduced in all ranks $j\geq i$, \item the sum of $t_i$-lengths 
of pieces of $q_3,$ between $<t_j>$-syllables for $j>i,$
is minimal for all paths $q_3$ with the above property.
\end{enumerate}
This number $s_1$ and the labels of the paths $q_3, q_4$ do not depend on $g.$
\end{cy}

The label of the path $q$ satisfying the conditions of this corollary will be 
called
a {\em $(u,u)$-pseudoconnector} for $h$, the label of the path $q_3$ will be 
called
a {\em $u$-pseudoconnector} for $h$ .  Notice that all elements of the
form $u^khu^s$ have the same $(u,u)$-pseudoconnectors and all the elements
of the form $hu^s$ have the same $u$-pseudoconnectors. 

\begin{df} In this definition we keep the notation of the corollary.
A {\em connecting zone} for  the element $h$  is defined
as follows.
Consider a path $q$ labelled by a $(u,u)$-pseudoconnector for $h$. 
Consider the  two 
phase vertices $b_1$ and $b_2$ closest 
to $p_2$ on either side of $p_2$ on the $u$-path $r_2,$  
with the property that the path from $p_2$ to $b_j$ has  $t_i$-length 
at least one. 
(If $p_2$ is not a phase vertex , then $b_1=gh\tilde u^{s_2}$,
$b_2=gh\tilde u^{s_2+k},$
where $k$ is either $1$ or $2$,
if $p_2$ is a phase vertex, 
then $b_1=gh\tilde u^{s_2-1}$
and $b_2$ is  $gh\tilde u^{s_2+1}$.) The connecting zone
for $h$ with respect to $g$ 
is the union of all phase vertices between such $b_1$ and $b_2$
for all $(u,u)$-pseudoconnectors ($u$-pseudoconnectors) for $h,$ see Fig. 3a .  
(Sometimes the connecting zone can consist of five phase vertices;
 they are shown in Fig. 3b.)

If the vertex $gh\tilde u^s$ belongs to the connecting zone for $h$ with respect 
to $g$, then the vertex $g_1h\tilde u^s$ belongs to the connecting zone for $h$ 
with respect 
to $g_1$. 
If it is clear from the context what the initial vertex $g$ of the
path
labelled by $h$ is (or if it does not matter), then 
we denote by $(h)_1$ and $(h)_2$ the initial
and terminal vertices of the connecting zone for $h$ with respect to $g,$
and will talk about the connecting zone for $h$ without mentioning 
the initial point $g$.

 \end{df}

\begin{center}{\bf 4. Construction of canonical representatives}\end{center}

We will define a section $\beta :H_{k_{i}+p}\rightarrow F_{k_i+p}$
(a mapping of sets such that $\pi\circ\beta =id$). For $X\in
H_{k_{i}+p}$ we call 
$\beta (X)$ the  {\em canonical representative of $X$.}

The canonical representative of an element in $G$ is just
some fixed geodesic word representing this element. 
(If $G$ is a free group, then 
it is the reduced word 
representing this element). 
Suppose we have already constructed representatives for all elements in 
$H_{k_i+p-1}.$
Now take an element $X$ in $H_{k_i+p}$
that contains $t_{k_i+p}$ 
and construct a representative
of this element in $F_{k_i+p}.$ 
The representative of $X$ will be a word (or a label of a 
path in $\Gamma (H_{k_i+p})$ ) 
corresponding to 
a reduced form of $X$ in  rank
$k_i+p.$ 
For each syllable $b$ of this element
between two consecutive $<t_{k_i+p}>$-syllables 
we will denote the canonical representative of the label of a path connecting
the vertices $(b^{-1})_2$ and $(b)_2$  
by $\bar b$ 
(see Fig. 4).  Suppose $X=b_0t^{\alpha _1}b_1t^{\alpha _2}\ldots b_n.$
Take some path labelled by $X$ in $\Gamma (H_{k_{i}+p})$. 
The path corresponding to 
$\bar b_i$ connects the points $(b_i^{-1})_2$
and $(b_i)_2$. Then the points
$(b_{i-1})_2$ and $(b_i^{-1})_2$ can be connected by a path
with the label $t_{i_k+p}^{\gamma _i}$. Then the canonical
representative for $X$ is $\bar b_0t_{i_k+p}^{\gamma _1}
\bar b_1t_{i_k+p}^{\gamma _2}\ldots \bar b_n$.

\begin{df} 
Let $\tilde u=u_{k_i+p}$, $h\in H_{k_i+p-1}$, $h$ is an element from
Corollary 2. The $(u,u)$-{\em connector} (resp. a $u$-connector) 
for the element $h$ is the canonical representative of the label of a
path, reduced in ranks $i,\ldots ,k_i+p-1,$ connecting the vertex 
$(h^{-1})_2$  (resp. $g$) with $(h)_2$ .

Let $\tilde u=u_j, \tilde v=u_k, j<k ,\ j,k\in\{{i+1},\ldots ,{k_{i+1}}\},$
and suppose that $q$ connects some vertex $g\tilde u^{s_1}$ with
$gh\tilde v^{s_2}, h\in H_i.$ Consider
the point $(h)_2$ (the terminal point of the path labelled by 
the $v$-connector for $h$). Let the path
$p$ connect $(h)_2$ with $g\tilde u^{s_1}$. Consider
the point $(\phi (p)^{-1})_2$ (the terminal point of the path labelled by 
the  $u$-connector for $\phi (p)^{-1}$).
The {\em $(u,v)$-connector} for $h$ is the canonical representative
of the label of a 
path in $H_i$ connecting
$(\phi(p)^{-1})_2$ with $(h)_2.$ 
\end{df}
It is important that 
the elements $u_j$ that do not belong to $G$
be chosen in such a way that
$\beta (u_{k_i+p}^n)=(\beta (u_{k_i+p}))^n.$  For this we first take
connectors for  $u_{k_i+p}^n$ and then, if necessary, replace
$u_{k_i+p}$ by its cyclic permutation  starting with $t_i.$
(If $G$ is a free group then
this equality can be made true for all $u_i$'s.)
Below, all the elements $u_j$'s will always be represented by the words 
$\beta (u_j)$ and we will write $u_j$ instead of 
$\beta (u_j)$. It will be clear from the context when $u_i$ means
a word and when it means the element represented by this word.

In Fig. 4 we considered the following example. Denote $t_{k_i+p}$
by $t$ and $u_{k_i+p}$ by $u$. Let $u=t^3$, $X=b_0t^{11}b_1t^{-7}b_2$
and $u=u_1u_2=u_3u_4=u_5u_6=u_7u_8$ in $H_i$. 
Let $b_0=c_0u_1^{-1}u^{-2}$, where $c_0$ is  a 
pseudoconnector for $b_0$,
 $b_1=uu_3c_1u_6u$, where $c_1$ is a pseudoconnector for $b_1$, 
$b_2=u^{-1}u_8^{-1}c_2$,
where $c_2$ is  a pseudoconnector for $b_2$. 
Let $\bar b_k$ ($k=0,1,2$)
be a connector corresponding to $b_k$. 
In this example, $\beta (X)=\beta (\bar b_0)t^8\beta (\bar b_1)t^{-10}
\beta (\bar b_2)$.

 \begin{center}{\bf 5. Middles}\end{center}

Let $X=\{X_1,\ldots ,X_L\}$ be a solution of system (\ref{1})
in the group $G^{{\bf Q}_{\pi}}.$ Suppose this solution contains
the minimal possible number of roots. 
Then for some $i$ this solution belongs to the group 
$H_{k_{i+1}}=H.$  Let ${\tau}_i (X)$ be the sum of the numbers of 
$t_i$-syllables in all $X_1,\ldots ,X_L.$ Denote 
$\tau (X)=({\tau}_1(X),\ldots ,{\tau}_{k_{i+1}}(X)).$
Let $X<Y$ if $\tau (X)<\tau (Y).$ 
Let $X=\{X_1,\ldots ,X_L\}$ be a minimal solution with respect to this
order.

By an {\em equational triangle} (resp. {\em equational diagram})
we mean a triangular equation $XYZ=1$ (resp. an equation $X_1\ldots X_n=1$) 
together with a solution $A,B,C$ (resp. $A_1,\ldots ,A_n$) 
and a diagram over $H$ having
$\beta (A)\beta (B)\beta (C)$ (resp. $\beta (A_1)\ldots \beta (A_n)$) 
as its boundary label. 
A {\em system of equational diagrams} is a system of equations
together with the system of diagrams, such that the solution associated to 
each equational diagram must be a solution of the whole system.

For the solution $X_1,\ldots ,X_L$ of system (\ref{1}) 
we will construct another solution
$X_1,'\ldots ,X_L'$ and a system of equations
over the group $G*K_{k_{i+1}}$ (if $G=F$, then $G*K_{k_{i+1}}=F_{k_{i+1}}$), 
such that $\beta (X_1'),\ldots ,\beta (X_L')$
will be a part of a solution of this new system, and every solution of the
new system will give a solution of system (\ref{1}).

Denote
$u_{k_{i}+p}$ by $u$ and $t_{k_{i}+p}$ by $t.$

Consider an equational triangle  in $H_{k_i+p}$ 
with at least one  side label containing $t$
(hence at least two side labels containing $t$). It is 
represented by a diagram in the form shown in Fig. 5 
(We showed only $t_j$-strips for $t_j=t$).

\begin{lm} \label{cr} 
Suppose that in $H_{k_{i+1}}$ we have a diagram (see Fig. 6)
with the boundary label $(b_1t_{j_1}^{r_{j_1}}\ldots t_{j_n}^{r_{j_n}}
b_nt_{j_{n+1}})
t_{j_{n+1}}^{r_{j_{n+1}}}(
c_1t_{j_1}^{p_{j_1}}\ldots t_{j_n}^{p_{j_n}}c_nt_{j_{n+1}})^{-1}$,
 where $j_1,\ldots ,j_{n+1}\in\{k_i+1,\ldots ,k_{i+1}\},$
 $b_1,c_1\ldots ,b_n,c_n\in H_{k_i}.$ 
 Then $\beta (b_1t_{j_1}^{r_{j_1}}\ldots t_{j_n}^{r_{j_n}}b_nt_{j_{n+1}})
t_{j_{n+1}}^{r_{j_{n+1}}}=
\beta (c_1t_{j_1}^{p_{j_1}}\ldots t_{j_n}^{p_{j_n}}c_nt_{j_{n+1}})$ in
$F_{k_{j+1}}$.
 \end{lm}
This assertion follows immediately from the choice of canonical representatives.

It follows from Lemma \ref{cr}
that every equational triangle either does not contain any cells or 
takes on one of the forms
shown in Fig. 7, and hence has a unique maximal nontrivial $H_{k_i+p-1}$-
subdiagram.

\begin{df} Consider an equational triangle in $\Gamma (H_{k_i+p})$ ($p\geq 1$) 
with the boundary label
$\beta (X_1)\beta (X_2)\beta (X_3).$
Suppose $X_1$ contains $t$.
Then the maximal nontrivial $H_{k_i+p-1}$-subdiagram of this triangle is called
the {\em middle} of the triangle. (The subdiagrams $ABCDEF$ (Fig. 7a), $ABCD$
 (Fig. 7b),
$ABC$ (Fig. 7c) are
{middles}.)

A boundary of a middle is canonically
subdivided into paths; each of the paths either belongs to
the $u$-side of a $t$-strip (is a $u$-path, joining two phase
vertices) 
or is a connector. 
The $u$-paths are called {\em pseudoangles} of the middle.
Every diagram over $H_{k_i+p}$ can be embedded into $\Gamma
(H_{k_i+p}).$ 
A {\em pseudoangle is long} if the corresponding strippath is nontrivial
and the connecting zones for two pseudoconnectors touching this pseudoangle do 
not
intersect, a {\em pseudoangle is short} 
if the corresponding strippath is nontrivial
and the zones of two pseudoconnectors touching this pseudoangle 
intersect, a {\em pseudoangle is trivial} if the
corresponding strippath is trivial.
A middle is called {\em triangular} if all the pseudoangles are trivial, it
is called short if it is not triangular and does not have long pseudoangles,
otherwise it is long.

If none of the $X_i$'s contains $t$ then the middle of the triangle coincides
with the triangle itself and is a {\em triangular middle}.
\end{df}
As an example, consider the middle $ABCDEF$ in Fig. 7a.
The paths $AB, CD, EF$  
are pseudoangles of this middle.

\begin{lm} \label{ABC} 
Consider an equational triangle over the group $H_{k_i+p}, p\geq 1,$ 
and the middle of this triangle. Let $AB$, $CD$ and $EF$ be the 
pseudoangles of the middle (see Fig. 7a) and $\phi (FA), \phi (BC),
\phi (DE)$ be the connectors. Suppose  $\phi
(AB)= u_{k_i+p}^n, n>0$. Then $\beta(\phi (FA)\phi (AB))=
 \phi (FA)\phi (AB)$. If $AB$ is a long pseudoangle and $B_1$ is
 the left end of the zone for $\phi (CB)$ then
 $A$ is the right end of the connecting zone for $\phi (EA)$,
 $B_1$ is the left end of the connecting zone for $\phi (DB)$ and 
 $$\beta(\phi (FA)\phi (AB_1)\phi (B_1C))=
 \phi (FA)\phi (AB_1)\beta\phi (B_1C),$$
  $$\beta(\phi (EA)\phi (AB_1)\phi (B_1C))=
 \beta \phi (EA)\phi (AB_1)\beta\phi (B_1C),$$
  $$\beta(\phi (FA)\phi (AB_1)\phi (B_1D))=
 \phi (FA)\phi (AB_1)\beta\phi (B_1D),$$
  $$\beta(\phi (EA)\phi (AB_1)\phi (B_1D))=
\beta \phi (EA)\phi (AB_1)\beta\phi (B_1D).$$ 
  \end{lm}
 
 {\bf Proof.}
 
 Let $HGG_1H_1$ be the paired $t_j$-strip closest to the point $A$ and 
$KLL_1K_1$
be the paired $t_{j_1}$-strip closest to the point $B$, $j,j_1>i$. If
there is no the strip $HGG_1H_1$, we just take $F=H_1, E=G_1$. If there
is no strip $KLL_1K_1$, we take $ K=D, L=C.$ 
We show this case in Fig. 8a. There is no $t_k$-strip , $k\geq i,$ connecting
the paths $H_1A$ and $BL$, because if there where, then from the description 
of representatives it would follow that $A=B.$ 

The label of the path $FH$ is the beginning of the canonical representative
for $\phi (FB)$. 
From the construction of representatives it follows that $H_1A$ contains
a $<t_i>$-syllable; let $PP_1$ be the $<t_i>$-syllable closest to the
vertex $A$. The path $AB$ begins with a $<t_i>$-syllable. Then 
$P_1A$ is the last $u_i$-connector in $u_{k_i+p}$. 
It follows from the choice of $u_{k_i+p}$ that $P_1A$ connects these 
$<t_i>$-syllables
and connects the $PP_1$-syllable with the vertex $A$ at the same time.
So $\phi (FA)\phi (AB)$ is the canonical representative of itself. 
Hence $\beta(\phi (FA)\phi (AB))=
 \phi (FA)\phi (AB)$.

Because the pseudoangle $AB$ is long,
there should be at least two paired $t_i$-strips
that begin on $AB$ and end on $ED$ as in Fig. 8b.
This implies that in the case when $D\not =C,$ 
the left end of the connecting zone for $\phi (DB)$ is $B_1,$
and in the case when $F\not =E$, 
the right end of the connecting zone for $\phi (EA)$ is $A.$

The four  equalities follow by symmetric considerations 
from the first one. The lemma is proved.

\begin{lm} Consider the diagram over the group $H_{k_i+p}, p\geq 1,$ 
shown in Fig. 7a or 7b. Suppose the middle of this diagram is not triangular.
\begin{enumerate}
\item If all three connectors contain some $t_j$'s greater than $t_i,$ then
two pseudoangles are trivial. 
\item If only two connectors contain some $t_j$ greater than $t_i$ then
one pseudoangle is trivial and another one is either trivial or short.
\item If none of the connectors contain
a $t_j$ greater than $t_i$ then two pseudoangles are either trivial
or short.
\end{enumerate}
\end{lm} 
{\bf Proof.} 
Consider the middle shown in Fig. 7a. Suppose there is a paired $t_j$-strip
$HGG_1H_1$ or a $t_j$-star $HGG_1I_1IH_1$ for $j>i,$ as shown in Fig. 8c. Let 
this 
strip be the strip
closest to $FE.$ If $j>k_i$ then by Lemma \ref{cr}
$E=F$ and this pseudoangle is trivial. If $j\leq k_i$, then the path
$HG$ does not contain any $t_k$ for $k\geq i.$  Let $E_1$ and $F_1$ be
the closest terminal points  of $(u_{k_i+p},u_{k_i+p})$-pseudoconnectors for 
$\phi (DE)$ and $\phi (AF)$ respectively.
 There are no paired $t_i$-arcs or $t_i$-stars
as shown by  broken lines in Fig. 8c. Hence the path $F_1E_1$
does not contain 
any $t_i$-syllable; hence $F=E$. If $HGG_1H_1$ is a paired $t_i$-strip,
then again the path $HG$ does not contain $t_i$, and this implies $F=E.$
If $HH_1G_1GJ_1J$ is a $t_i$-star, then $F_1$ and $E_1$ must coincide,
because $AE_1$ is also a $(u_{k_i+p},u_{k_i+p})$-pseudoconnector for $AF$. 
Now if all the three connectors contain some $t_j$'s greater than $t_i,$ then
either there are  at least two paired strips for $t_j$
and for $t_{j_1}$, as shown in Fig. 8a,  or there is a $t_j$-star. Then $F=E$, 
$D=C$. 
So  at least two pseudoangles are trivial. 

If only two connectors contain some $t_j$ greater than $t_i,$ then
there is a paired $t_j$-strip $HGG_1H_1$ closest to the vertex $A$ 
as in Fig. 8c; hence $F=E$ and this
pseudoangle is trivial. Suppose that the pseudoangle $CD$ is nontrivial,
the proof that it is short is similar to the proof of the third
assertion that follows below. 

If none of the connectors contains a $t_j$ greater that $t_i$, then
there are no paired $t_j$-strips in the middle, and the $t_i$-strips
are situated as in Fig. 9 a,b,c,d,e. In Fig. 9a the pseudoangles
$FE$ and $DC$ are trivial,
in Fig. 9b the pseudoangles
$FE$ $AB$ and $DC$ are short, in Fig. 9c the pseudoangle
$FE$ is short and $DC$ is trivial, in Fig.~9d the pseudoangle $FE$ is short,
$CD$ is trivial. 
In the case shown in Fig. 9e two pseudoangles are trivial and one is
short. 
Indeed, $A_1$ and $A_2$ are both the
beginnings of some $(u_{k_i+p},u_{k_i+p})$-pseudoconnectors for $BC$,
$D_1$ and $C_1$ are both the beginnings of some
$(u_{k_i+p},u_{k_i+p})$-pseudoconnectors for $CB$.

The lemma is proved.

\begin{center}{\bf 6. Reduction to the free group}\end{center}
We start with a system $S$ of $M$ equational triangles in rank $k_{i+1}.$
We will construct a system of equations in $G*K_{k_{i+1}}$ such that
the solution of $S$ is a part of a solution of this system.

Consider in $S$ an equational triangle with the boundary $X_1X_2X_3$ 
corresponding to Figure 10. Recall that $\beta (X_i)=X_i$ in $F_{k_{i+1}}$. 
Then
$$X_1=P_1Y_1(t^{s_{k_i+1}m_1+r_1})P_2,$$ 
$$X_2=P_2^{-1}t^{s_{k_i+1} m_2+r_2}Y_2^{-1}P_3^{-1},$$
$$X_3=P_3Y_3P_1^{-1}$$ 
$$Y_3=Y_4u^{-(m_1+m_2)}Y_1^{-1}$$
in the group $F_{k_{i+1}},$
and we have to add the conditions that $u\in F_i$ , $r_1+r_2=s_i
$, as well as the 
equational triangle over $H_{k_{i+1}-1}$ with the boundary $Y_4uY_2^{-1}.$   

\begin{df} A {\em free equational diagram} (resp. {\em triangle}) 
is an equational diagram (resp. triangle)
with no cells. \end{df}
The 
equational triangle with the boundary  $Y_4uY_2^{-1}$ over $H_{k_{i+1}-1}$  
is ``non-free'' (it means that the corresponding diagram might contain
cells). We will say that the first three free equational diagrams are
{\em free equational diagrams (or free equations ) of the first type} and the 
fourth free 
equational  diagram is a {\em free equational diagram (or free equation)  
of the second type}.
We also can rewrite free equational triangles in $F_{k_{i+1}}$ in the form
$$X_1=P_1Y_1T^{ m_1}T_1P_2,$$ 
$$X_2=P_2^{-1}T^{ m_2}T_2(Y_2)^{-1}P_3^{-1},$$
$$X_3=P_3Y_3P_1^{-1}$$ 
$$Y_3=Y_4u^{-(m_1+m_2)}Y_1^{-1}$$
in the group $F_{k_{i+1}},$
and add equations $$[T,t]=[ T_1,t]=1,[T_2,t]=1,\ \ T_1T_2=T_1.$$
We  also have to add the condition that $u\in F_i$ .

Notice that the new ``variables'' in the non-free equational triangles  
are $u,$  the connector $Y_2,$ participating in the middle, 
and the canonical representative of the 
label of the path connecting $C$ with the left end of the
zone for $\phi (CB)$ ($\beta\phi (Y_2u^{-1})$). Notice also that 
the freeness of the equational triangles of the first type imply that some of 
the 
$Y's$
are the same. Below we will often use the fact that they will still have a 
solution in $X's$ and $P's$  
if we replace some of the $Y's$ by some other words.

We can rewrite every equational triangle in a similar way.
We will obtain
a system of free equations over $F_{k_{i+1}}$ and a system of
triangular equations
over $H_{k_{i+1}-1}$. Denote the union of these systems by $S_1.$
At the next step we repeat this process for the non-free equations
over $H_{k_{i+1}-1}$ with unknowns $Y_i$. Finally, after $k_{i+1}$ steps
we will obtain a system $S_{k_{i+1}}$ over a group $G*K_{k_{i+1}}$,
which, in the case when $G$ is free, is just a 
free group $F_{k_{i+1}}$. Every
solution of $S_{k_{i+1}}$ gives a solution of the system (\ref{1}).

\begin{df} For a given system of equational triangles in $H$, the union of 
all non-free equational triangles, in all
ranks from $0$ to $k_{i+1}$ 
will be called a {\em tower of non-free equational triangles 
generated by this system.}
The same non-free equational triangle can occur in different ranks; 
in this case we consider all
these occurrences as one triangle of the tower. \end{df} 
\begin{lm} If a solution of a system of triangular equations is minimal,
then in the tower 
of non-free equational triangles generated by this system
there is at least one nontriangular middle
for the non-free triangles  in each rank . \end{lm}
{\bf Proof.} Suppose that in rank $j$ all the middles are triangular. 
Let $\beta (X_1),\ldots ,\beta (X_L)$ be the canonical
representatives of the elements of the solution. Then they are part
of a solution of the system $S_{k_{i+1}}$. 
The constant $t_i$ does not participate in the equations
of the system $S_{k_{i+1}}$. Then we can cut out $<t_j>$-syllables 
from all the elements in the solution of $S_{k_{i+1}}$
and again have a solution of $S_{k_{i+1}}$ .
Every solution of $S_{k_{i+1}}$ produces
a solution of the system (\ref{1}).
This contradicts the minimality of the original solution.
\begin{lm} \label{one}
Every  middle in rank $k_i+p$ can be split
into several free equational triangles and at most one non-free equational
triangle in the previous rank.\end{lm}

{\bf Proof.} 
Every triangular middle will be itself a triangular equation in the
previous rank.

Consider the equational triangle in Fig. 7a . Let ABCDEF be the middle of 
this  triangle. Suppose the pseudoangle
$AB$ is nontrivial and $AB$ is a positive power of $u$. 
Let $B_1$ be the left end of the zone for $\phi (CB)$,
and $B_2$ be the left end of the zone for $\phi (DB)$.
Consider the following possibilities. They all are shown in Fig. 11
with the corresponding letter; in the cases a)-d) the pseudoangle
$AB$ is long. 
Non-free equational triangles are marked by a 
star. All the assertions below follow directly from Lemma \ref{ABC}.

a) Two other pseudoangles are trivial . We have $E=F, C=D.$
Then there is a free equational diagram, with the boundary
$\phi (FA)\phi (AB_1)\beta\phi (B_1C)\phi (CF),$ and an equational
triangle with the boundary $\phi (CB)\phi (BB_1)\beta\phi (B_1C).$

b) The pseudoangle $EF$ is trivial, $\phi (CD)$ is a negative power of $u$.
There is one non-free equational triangle, 
with the boundary $\phi (CB)\phi (BB_1)\beta\phi (B_1C),$ and two free 
equational 
diagrams with the boundaries $\phi (CD)\phi (DF)\phi (FD),$ and 
$\phi (FA)\phi (AB_1)\beta\phi (B_1C).$

c) The pseudoangle $CD$ is trivial, $\phi (FE)$ is a negative power of $u$.
There is one non-free equational triangle, 
with the boundary $\phi (CB)\phi (BB_1)\beta\phi (B_1C),$ and two free 
equational 
diagrams with the boundaries $\phi (CE)\phi (EF)\phi (FC)$ and 
$\phi (FA)\phi (AB_1)\beta\phi (B_1C).$
 
d) The pseudoangles $CD$ and $FE$ are both nontrivial, $CD$ and $FE$ are
negative powers of $u$.
There is one non-free equational triangle, 
with the boundary $\phi (CB)\phi (BB_1)\beta\phi (B_1C),$ and two free 
equational 
diagrams with the boundaries $\phi (CD)\phi (DE)\phi (EF)\phi (FC)$ and 
$\phi (FA)\phi (AB_1)\beta\phi (B_1C).$
 
In the case where there are no long pseudoangles, the middle is short,
and, up to the relabelling of the vertices, there are two different 
possibilities:

e) $EF$ and $DC$ are positive powers of $u$.
There is one non-free equational triangle, 
with the boundary $\phi (CB)\beta \phi (BF)\beta\phi (FC),$ and two free  
equational 
diagrams with the boundaries  
$\phi (EF)\phi (FC)\phi (CD)\phi (DE)$ and $\beta\phi (FB)\phi (BA)\phi (AF).$

f) $EF$ and $CD$ are positive powers of $u$.
There is one non-free equational triangle, 
with the boundary $\beta\phi (FB)\beta\phi (BD)\beta\phi (DF),$ and three free 
equational 
diagrams with the boundaries $\phi (CD)\beta\phi (DB)\phi (BC)$,
$\phi (EF)\beta\phi (FD)\phi (DE)$  and  
$\phi (FA)\phi (AB)\beta\phi (BF).$

The lemma is proved. 

\begin{lm}\label{steps} Consider $n$ non-free equational triangles
in rank $k_i+p$. Suppose the solution is minimal. 
 Then we can rewrite the system as a system of free 
equational triangles and at most $n$  non-free equational triangles
in ranks less than $k_i+p$.
The variables in the new non-free equational triangles are of the
following types: the connectors in the middles in rank $k_i+p$,
and variables $\bar Z_q$,  
such that 
there is a free equation, 
either of the form $\bar Z_q=u_{k_i+p}^q$  
or of the form $\bar Z_q=u^{q_1}Tu^{q_2}$, where $T$ is the label of some
connector.  
\end{lm}
{\bf Proof.}
Every nontriangular middle can be split by the previous lemma
into several free equational triangles and at most one non-free equational
triangle. The form of the new variables is obtained in the proof of 
the previous lemma.

\begin{lm}\label{verygood} 
In the minimal solution, the difference $k_{j+1}-j,$ 
and hence the number of $u$'s containing some fixed $t_j$ as the greatest root
(which is equal to the difference $k_{j+1}-k_j$) cannot
be more then $3M,$ where $M$ is the number of equations in the original
system. 
\end{lm}

{\bf Proof.}
 We will show that, 
 when we consider
 non-free equational triangles in ranks $k_{j+1},k_{j+1}-1,\ldots $,
   we have to come to rank $j$  after not more than $3M$ steps.
Indeed,  
consider one  non-free equational triangle in rank $k_{j+1}$. In Fig.
12 we show three paired $t_p$-strips, for $p>i,$ closest to the center
of this triangle.  Let them be $t_{p_1},t_{p_2}$ and $t_{p_3}$-strip.
Then this triangle will produce non-free triangular equations with
nontriangular middles  in  at most three ranks from $k_{j+1}$ to $j+1$.
These ranks can be only $p_1,p_2$ and $p_3$. 
In each rank  we have at 
least one equational triangle with non-triangular middle.
Hence the maximal possible number of ranks between $k_{j+1}$ and
$j+1$ is $3M$. 
Lemma is proved.

After considering some level $j,$ we have an effectively bounded number of 
free equations of the second type
connecting the variables $Z_1,\ldots , Z_m,$  formed on
this level, and powers of $u_k$'s, where $k\in\{k_j,\ldots ,k_{j+1}-1\}.$
These equations have one of the following form
$Z_{j_1}=Z_{j_2}u_k^qZ_{j_3}$ or 
$Z_j=u_k^l, Z_j=u_k^lZ_su_k^r,$ where $l,r\in\{1,2,3,4\}$. 

\begin{center}{\bf 7. Shortening of middle-strips}
\end{center}
We now  want to bound the powers of $u_j$'s in the free equational diagrams.
Here $j\geq k_1,$ in the case when $G$ is not a free group, and $j\geq 1$
in the case when $G$ is a free group.

We denote $u=u_{k_{i+1}}, t=t_{k_{i+1}}.$

\begin{df} A {middle-$t$-strip} is a $t$-strip formed by $t$-cells 
such that
their $u$- or $t$-sides belong to the 
long pseudoangle of the middle of 
some non-free equational triangle in rank $k_{i+1}.$
\end{df}

Notice that if two paths contain $t_i,$ and 
represent the same
element in $H_{i},$ 
then their canonical representatives are chosen in such
a way that the part between the first and the last $<t_i>$-syllable
is the same word in $F_{k_{i+1}}$.
So we can talk about the subword $u$ in $H_{i}.$ (Recall that we
write $u$ instead of $\beta (u).$)

\begin{df} A subword $u^k$  
of the label of a side of an equational triangle is called {\em shrinkable} if
one of the following conditions is satisfied
\begin{enumerate}
\item
There is an occurrence of $u$ in $u^k$ which is on the boundary
of some $t$-cell belonging to some middle-$t$-strip.
\item $u^k$ is a subword of some shrinkable subword $u^l$ of a side
of an equational triangle.
\item Suppose there is a common subpath of two sides of an equational triangle
labelled by $\beta (X)$ and $\beta (Y)$. Suppose that the label of this common 
subpath is a shrinkable
subword of $\beta (X)$, then it is also a shrinkable subword of $\beta (Y).$
\end{enumerate}
So, to be precise, the set of all shrinkable subwords  of labels $\beta (X)$ of 
sides of 
equational triangles is the smallest class of subwords satisfying the
above description.
\end{df}

A similar definition can be given for shrinkable subwords $t^k$.

Suppose we have an equational triangle with the boundary label $XYZ$. Suppose 
this
triangle contains a middle-$t$-strip. Then some maximal shrinkable subword 
$u^k$ of some variable, say $X$, 
contains a part of the boundary of this middle-$t$-strip. There are
two possibilities. Either some
maximal shrinkable subwords $t^{k_{i+1}n_1}$ of $Y$ and $t^{k_{i+1}n_2}$
of $Z$   both contain a part of the boundary of this middle-$t$-strip
(Fig. 13a) or the maximal shrinkable subword $t^{k_{i+1}n_1}$ of $Y$
contains the $t$-side of the boundary of this middle-$t$-strip, and
the maximal shrinkable subword $t^{k_{i+1}n_2}$ of $Z$ is the label of a
common path of $Y$ and $Z$ (going along
the other part of the subword $t^{k_{i+1}n_1}$ of $Y$) (Fig. 13b).

Suppose the triangle $XYZ$
does not contain a middle-$t$-strip. It might happen that one part of 
some maximal
shrinkable subword $u^k$ (resp. $t^{ks_{k_{i+1}}}$) of $X$ labels
some common path of $X$ and $Y$  and another part labels some
common path of $X$ and $Z$
(see Fig. 13c,d). 
In both  cases {\em a maximal shrinkable piece} of the triangle 
is the minimal subdiagram  having these three
maximal shrinkable subwords on its boundary. (In the cases shown
in Fig. 13c,d, these pieces $ABCD$ and $ABC$ do not contain any cells.) 

Denote the set of maximal shrinkable pieces by $\cal S.$ 
The following lemma is obvious.
\begin{lm}
Every maximal shrinkable piece has one of the forms shown in Fig. 14 (broken 
line 
for $t$,
bold line for $u$).
\end{lm}

Below, we will consider only shrinkable subwords $u^k$, $k\geq 4$ and
$t^l$, $l> 4s_{k_{i+1}}.$
\begin{df} The {\em length of the shrinkable subword $u^k$} is $k-4$, the
{\em length of the shrinkable subword $t^{s_{k_{i+1}}j+r}$,} 
where $0<r<s_{k_{i+1}},$ is $j-4$. 
\end{df}
  
Every maximal shrinkable piece contains three maximal shrinkable subwords. 
Denote them by $a_{kj}$, 
$j=1,2,3,$ and their lengths
by $\ell _{i1},
\ell _{i2},\ell _{i3}.$ For example, $a_{k3}=u^{4+\ell _{k3}},
a_{k1}=t^{s_{k_{i+1}}(4+\ell _{k1})}, a_{k2}=t^{s_{k_{i+1}}(4+\ell _{k2})}.$
Each diagram in $\cal S$ gives one of the equations
$\ell _{i2}+s=
\ell _{i3}$ 
, where $s\in\{4,5,6,7,8\}$. Denote this system of equations, assigned 
to $\cal S,$ by
$\cal L.$ Add to this system (with less than 3M unknowns) the equalities
and inequalities $\ell _{ij}=\ell _{tp},$ $\ell _{ij}>0, \ell _{ij}=0,$
and the equations that keep the length of the short pseudoangles
and maximal shrinkable subwords of nonpositive length.

It is possible to obtain only a finite number (depending on $M$) of
distinct linear systems of this form. Every such system is 
algorithmically solvable. This follows from the fundamental result of
Presburger \cite{P} about the decidability of the elementary theory
of the natural numbers with addition.
Take one 
solution for each  system. Let $\bar L$ be the maximum of the $\ell_{ij}$
in these solutions. Then the system $\cal L$ also has a solution
bounded by $\bar L.$ Let
$$\{\bar\ell _{kj}, i=1,\ldots ,M, j=1,2,3\} $$
be this solution, $b_{kj}$ the pieces corresponding to $a_{kj},$ but of length
$\bar\ell _{kj}$ (for example, if $a_{k3}=u^{4+l_{k3}}$, then
$b_{k3}=u^{4+\bar l_{k3}}$.) 

\begin{lm} 
\label{del} Replace all the pieces $a_{ij}$ by the analogous pieces
$b_{ij}.$ We get another solution of the original system of equations.
\end{lm}
Suppose we replaced all $a_{ij}$ by the corresponding $b_{ij}$. If
some $a_{ij}$ participated in several maximal shrinkable pieces , 
then $l_{ij},$ and 
hence $\bar l_{ij},$ satisfy linear equations for all these pieces.
In the places which do not meet shrinkable pieces, replacing $a_{ij}$
by $b_{ij}$ also does not destroy the solution. The lemma is proved.

We described how to obtain another solution with the length of middle
strips in rank $k_{i+1}$ bounded by $\bar L+4$. We now use induction. 
Suppose that the middle strips are bounded for all 
$u_{k_{i+1}-1}\ldots u_{k_{i+1}-j+1},$ and bound them for
$u_{k_{i+1}-j}$.
We add new variables and rewrite free equational diagrams of the second type as
equational triangles. 
We will obtain a bounded number of equational triangles,
because the length of the middle strips for 
$u_{k_{i+1}-1}\ldots u_{k_{i+1}-j+1}$ is bounded and the number of 
equational triangles obtained from free equational diagrams 
depends on the lengths
of middle strips.

We got some shrinkable pieces in rank $k_{i+1}-j$; 
we name then $a_{rl},$ $l=1,2,3.$ 
If $a_{rl}$ coincides with a piece of another solution not in the middle
we also name this piece $a_{rl}.$
The word $a_{rl}=u_{k_{i+1}-j}^m$ can be a subword of some
greater word $u_{k_{i+1}-p}$ 
and touch it somewhere, but since we consider all free and non-free equational
triangles,
this is taken into consideration in the equations.
The replacement of $a_{rl}$ in $u_{k_{i+1}-p}$ by $b_{rl}$
can turn $u_{k_{i+1}-p}$ into a proper power (in which case we turn 
the corresponding
$t_{k_{i+1}-p}$ into the same power ), 
but it cannot turn $u_{k_{i+1}-p}$ and $u_{k_{i+1}-{p_1}}$ into
proper powers of the same word, because the number of roots was taken to be 
minimal.

Indeed, if we obtained $u=v^s, u=t_1^{\alpha}, v=t_2^{\gamma},$
we could put $t_1=t_3^{s{\gamma}}, t_2=t_3^{\alpha}$ and this would be another
 solution with  fewer   roots added.

Now suppose by induction that we bounded the lengths of the middle-strips at
all levels
higher than $j.$ 
Free equational triangles at levels higher that $j$ do not
affect the maximal shrinkable pieces at level $j$. Indeed, free equations on 
levels higher than $j$ just equate some of the
$<t_r>$-syllables, for $r>k_{j+1},$  and some
of the $(u_k,u_l)$-connectors between two neighboring $<t_k>,<t_l>$-syllables,
where $k,l>j$ are equal. All the maximal shrinkable subwords $u_j^k$
 at level $j$ are 
inside these connectors. All the maximal shrinkable subwords
$t_r^k$, where $r\in\{k_{j}+1,\ldots ,k_{j+1}\},$ coincide with some of these
$<t_r>$-syllables.
For the level $j$  we just repeat the procedure described for 
the level $i.$
 Because there is only a bounded number of $u$'s on the same level,
and the number of non-free equations does not exceed $M,$ 
we can get a boundary for the lengths of middles at each level.
So all middle $t$-strips, at all the levels greater than $0,$ can be shortened .

If $G$ is a free 
group, the
middle $t$-strips at level $0$ can be shortened by the same procedure.

\begin{lm} \label{deep} 
Starting with a minimal solution, with  
bounded lengths of middle $t_k$-strips at all levels, 
it is possible to  construct a solution with bounded depth in every
rank $j$ 
of the root $t_j$ 
(i.e. the $s_j$'s are bounded).
\end{lm}
{\bf Proof.}  
Suppose, by induction, that we have bounded the depth of the roots on 
all levels
higher than $j,$ and in the ranks $k_{j+1},\ldots , k_{j+1}-p+1$ .
Now we will bound it in rank $k_{j+1}-p.$ Consider all non-free
equational triangles in rank $k_{j+1}-p,$ together with  all their middle strips
as shown in Fig. 14a,b,c. Let $a_{ik}=t_{k_{j+i}-p}^{4+l_{ik}}$ be a shrinkable
$t_{k_{j+i}-p}$-subword. 
For the middle strips we can write the
system of linear equations with variables $l_{ik}$ and $d$, where $d$
is the depth of the root $t_{k_{j+1}-p}.$

We have no more than $M$ equations of the form 
$t_{k_{j+1}-p}^{x_1}t_{k_{j+1}-p}^{x_2}=u_{k_{j+1}-p}^s$ ($s$ is bounded by the
maximal length of middle-strips)
 or $t_{k_{j+i}-p}^{x_1}t_{k_{j+i}-p}^{x_2}=t_{k_{j+i}-p}^{x_3}$. So we have 
a homogeneous system of linear equations with not more than $3M$ variables . 
These 
equations can only be of the form $x_i\pm x_j=sd$ 
($s$ is bounded, by $\bar L$), or $x_i+x_j=x_k$ 
The prime divisors of $d$ must belong to the set $\pi$, $x_i\geq 0$. 
The minimality 
of the solution $X=\{X_1,\ldots ,X_L\}$ implies that the linear system does 
not admit solutions in which one of the $x's$ is divisible by $d$. 
This implies 
that, after transforming the system to reduced row-echelon form, we can only 
have $d$ as a free variable, and can either solve the system  or see that it is 
unsolvable.

There can be only a finite number of such systems, with each system having some 
solution;
if $\tilde L$ is the greatest value of $d$ in all these solutions, then
we always can find a solution of system (\ref{1}) where the depth
of the roots is bounded by $\tilde L.$
The lemma is proved.
\begin{center}{\bf 8. The proof of Theorem 1}
\end{center}

From now on we suppose that the group $G$ is free. So $F=G,$ and
on level zero  we have only free equations. 
\begin{prop}
\label{n}
It is possible to determine a recursive function $\psi (M)$ such that 
if a system of $M$ triangular equations has a solution in $F^{\bf Q}$, then
it also has a solution in some group $H_{\psi (n)},$ where $H_0=F$ and
$H_{i+1}$ is obtained from $H_i$ by adding some root.
\end{prop}  

{\bf Proof.} 
It is possible to find a number $\psi (M)$ such that if $k_{i+1}$
 is greater than $\psi (M)$, there will be some numbers $r$ and $s<r$ such 
 that $k_{r+1}-r=k_{s+1}-s$ and 
the following conditions are satisfied: \begin{enumerate}
\item The systems of non-free equations in ranks $r$ and $s$ are equivalent
(there exists a bijection between their sets of variables that
induces a bijection between their solution sets). 
\item 
Let $Y_1,\ldots ,Y_p,u_{r+1},\ldots ,u_{k_{r+1}}$ be the variables
in the system of non-free equations in rank $r$ 
and $Z_1,\ldots ,Z_p,u_{s+1},\ldots ,u_{k_{s+1}}$ be the variables
in the system of non-free equations in rank $s$.
As variables,  $u_{r+1},\ldots ,u_{k_{r+1}}$ correspond to
$u_{s+1},\ldots ,u_{k_{s+1}}$.
\item The  corresponding 
roots of  $u_{r+1},\ldots ,u_{k_{r+1}}$ and
of $u_{s+1},\ldots ,u_{k_{s+1}}$ have the same depth.
\item Consider first free equations of the second type in ranks
from ${r+1}$ to $k_{r+1}.$ They imply certain free equations for
$(u_j,u_k)$- and $u_j$-connectors in $H_r$
(those which are
defined) and 
original variables in $H_r$ 
, $j,k\in
\{r+1,\ldots ,k_{r+1}\}$. Hence the variables $Y_1,\ldots
,Y_p,u_{r+1},\ldots ,u_{k_{r+1}}$ satisfy some free equations. 
Next, consider the free
equations of the second type in ranks
from ${s+1}$ to $k_{s+1}.$ They likewise imply certain free equations
for
$(u_j,u_k)-$ and $u_j$-connectors and 
original variables in $H_s$
, $j,k\in
\{s+1,\ldots ,k_{s+1}\}$. Hence the variables $Z_1,\ldots
,Z_p,u_{s+1},\ldots ,u_{k_{s+1}}$ also satisfy  some free equations. 
The systems of the free equations for $Y_1,\ldots
,Y_p,u_{r+1},\ldots ,u_{k_{r+1}}$ and for $Z_1,\ldots
,Z_p,u_{s+1},\ldots ,u_{k_{s+1}}$ must be equivalent.
 \end{enumerate}
Indeed, the number of roots on a given level, as well as the difference
$k_{j+1}-j,$ are bounded by Lemma
\ref{verygood}, 
the length of middlestrips is bounded, the number
of free equations of the second type at every level is bounded,
and the depth of the roots is bounded by Lemma \ref{deep}.

We replace now all the variables in equational triangles in rank $r$ 
by the corresponding variables
in equational triangles in rank $s.$ 
In particular $u_{r+1},\ldots ,u_{k_{r+1}}$ will be replaced
by $u_{s+1},\ldots ,u_{k_{s+1}};$ $(u_k,u_l)-$ and $u_k$-connectors,  for 
$k,l\in
\{r+1,\ldots , k_{r+1}\}$ participating in non-free equational triangles,
will be
replaced by the corresponding $(u_k,u_l)-$ and $u_k$-connectors,  for $k,l\in
\{s+1,\ldots , k_{s+1}\}.$ If some 
$(u_k,u_l)$ or $u_k$-connector,  for $k,l\in
\{r+1,\ldots , k_{r+1}\},$ or original variable in $H_r$ 
does not participate in nonfree equations, then
we can replace it by an arbitrary element in $H_s.$

Our purpose is to show that we will  again obtain a solution.
Indeed, consider first what will happen to the free equational triangles  
that come from the levels higher than $r.$ These equations just
indicate that some of the $(u_k,u_j)$- and $u_j$-connectors 
are the same and some of the $<t_j>$-
syllables are the same,
for
$j,k\in\{r+1,\ldots ,k_{i+1}\}$ . We will have to replace
the pieces, corresponding to 
$u_{r+1},\ldots ,u_{k_{r+1}}$ in these equations, by the pieces
corresponging to $u_{s+1},\ldots ,u_{k_{s+1}}$ , 
and hence to change $u_{k_{r+1}+1}$ and all the highest 
$u$'s. But all these pieces are between two neighboring occurrences of
$t_j$ and $t_k$, $j,k\in\{r+1,\ldots ,k_{i+1}\}$, 
so do not affect the other structure of free equations.
We also will replace the pieces corresponding
to $(u_j,u_k)$-connectors in $H_r,$ for $j,k\in\{r+1,\ldots ,k_{r+1}\},$ by
the pieces corresponding
to $(u_j,u_k)$-connectors in $H_s,$ for $j,k\in\{s+1,\ldots ,k_{s+1}\},$  
and again this will change some connectors;however, all the connectors which
were the same will remain the same.

So we obtain a solution with fewer roots (we do not need 
$t_{s+1},\ldots , t_{k_{s+1}}$, if $k_{s+1} <{r+1}$ or 
$t_{s+1},\ldots , t_{{r+1}-1}$ if $k_{s+1} \geq {r+1}$). 
Indeed, now $t_{r+1},\ldots , t_{k_{r+1}}$ will now be the roots
of $u_{s+1},\ldots , u_{k_{s+1}}$.
We obtained a contradiction with the 
minimality of our solution.  

The proposition is proved.

It follows from the Proposition that $k_{i+1}=\psi (M)$. The length
of the middle $t_j$-strips, in all the ranks $j,$ is bounded. We will construct 
a
system of equational triangles in $F_{k_{i+1}},$ as 
it was described in section 6.
The  number of possible systems of free equations that correspond to
these systems of equational triangles  is bounded, and we can list
them.
Finally we have a finite number of possible systems of 
equations  in the free group $F_{k_{i+1}},$
with the restriction that some of the variables belong to subgroups
generated by only a part of the generating set of the free
group $F_{k_{i+1}}$. 
Each such system is algorithmically decidable \cite{Mak84}.
If all these systems are
incompatible, then our system does not have a solution. If at least
one of them is consistent, we obtain a solution of our system by
substituting the corresponding $P_i$'s $u_j$'s and $T_k$'s in the expressions
for unknowns.
\begin{center}{\bf 9. The proof of Theorem 2}
\end{center}

To prove Theorem 2 we first need some definitions.

A subgroup $G$ is called {\it existentially closed } in a group $H$
if any existential sentence, with constants from $G,$ holds in the whole
group $H$ if and only if it holds in the subgroup $G$.

Let $G$ be a subgroup of $H$. A finite system of equations
 $W(x,g) = \{ w_1(x,g) = 1, \ldots , w_k(x,g) = 1 \},$ 
with variables $x = (x_1, \ldots, x_n)$ and constants $g = (g_1, \ldots , g_m)$ 
from $G$, has a solution in $H$ (resp. in $G$) if and only if the following 
formulae holds
in $H$ (resp. in $G$):

$$\exists x ( w_1(x,g) = 1 \wedge \ldots  \wedge w_k(x,g) = 1).$$

Therefore, If $G$ is existentially closed in $H$, then any system $W(x,g) = 1$, 
with 
constants in $G$, has a solution in $H$ if and only if it has a solution in $G$.

Let $G$ be a subgroup of $H$. Following \cite{BMR} we will say that  {\em  $H$ 
is locally $\omega$-separated in $G$ by retractions } if for arbitrary finitely 
many 
  nontrivial elements $h_1, \ldots, h_n \in H$ there exists a homomorphism
$\psi: H \longrightarrow G$, which is  the identity on $G$, such that the images 
of 
$h_1, \ldots, h_n $ under $\psi$ are also nontrivial in $G$.

In \cite{BMR} the following result has been proven:
 Let $G$
be a torsion-free hyperbolic group and $A$ a ring of characteristic 0. 
 Then $G^A$ is locally $\omega$-separated in $G^{{\bf Q}_{\pi(A)}}$ by 
retractions.

\begin{lm}
Let $G$ be a subgroup of $H$. If $H$ 
is locally $\omega$-separated in $G$ by retractions, then 
$G$ is  existentially closed  in  $H$.
\end{lm}
{\bf Proof.} From general predicate calculus we know that any existential 
sentence 
in the group theory language
with constants from $G$ is equivalent to a sentence of the following type:

$$ \Phi = \exists x(\bigwedge_{1}^{s}u_i(x,g) = 1 \bigwedge_{1}^{t}v_j(x,g) \neq 
1),$$
where the $u_i$'s and  $v_j$'s are group words, $x = (x_1, \ldots, x_n)$ are  
variables and 
$g = (g_1, \ldots , g_k)$ are some constants from $G$. 

Let the elements $h = (h_1, \ldots, h_n) \in H$ satisfy the quantifier-free part 
of this
 sentence in $H$.
Denote by $H_0$ the subgroup $<h_1, \ldots, h_n> $ in $H$.  By the conditions
of the lemma there exists a homomorphism $f :  H_0 \longrightarrow G,$ 
which separates the elements $v_1(h_1, \ldots, h_n), \ldots,
v_t(h_1, \ldots, h_n)$ in $G$.  This implies that the images $f(h_1),
\ldots, f(h_n)$ satisfy in $G$ the same equalities $u_i(f(h_1), \ldots,
f(h_n)) = 1, i = 1, \ldots, s,$ and inequalities $v_i(f(h_1), \ldots,
f(h_n)) \neq 1, i = 1, \ldots, t.$  Therefore, the sentence
$\Phi$ holds in $G$.  This shows that $G$ is existentially-closed in $H$. 
$\Box$

\begin{cy}
Let $G$
be a torsion-free hyperbolic group and $A$ a ring of characteristic 0. 
 Then $G^{{\bf Q}_{\pi(A)}}$ is existentially closed in  $G^A$ .
\end{cy}

Now we can complete the proof of Theorem 2. 
By the corollary above, $F^{{\bf Q}_{\pi(A)}}$ is existentially 
closed in $F^A$, hence any system $W = 1$ with coefficients from $F$ has a 
solution in 
$F^A$ if and only if it has a solution in $F^{{\bf Q}_{\pi(A)}}$. The result 
now follows from 
Theorem 1.  

\begin{center}{\bf 10. Equations and inequalities in $F^{{\bf Q}_{\pi}}$}
\end{center}
In this section we show how to reduce the problem about algorithmic 
decidability of the universal theory of the group $F^{{\bf Q}_{\pi}}$
(with constants from $F^{{\bf Q}_{\pi}}$) to some question in the free group 
$F$.

To prove the decidability of the universal theory of $F^{{\bf Q}_{\pi}}$
one has to construct the algorithm solving finite systems of equations and
inequalities in $F^{{\bf Q}_{\pi}}$. Such a system can be reduced to
the system of triangular equations (\ref{1}) together with
some inequalities $X_k\not = 1.$ The algorithm constructed in this paper
can be used to reduce the system to a finite number of possible systems of 
equations  and inequalities in the free group $F_{k_{i+1}},$
with the restriction that some of the variables belong to subgroups
generated by only a part of the generating set of the free
group $F_{k_{i+1}}$.


\end{document}